\documentclass[12pt, reqno]{amsart}
\usepackage{amsmath, amsthm, amscd, amsfonts, amssymb, graphicx, color}
\usepackage[bookmarksnumbered, colorlinks, plainpages]{hyperref}

\newtheorem{theorem}{Theorem}[section]

\newtheorem{proposition}[theorem]{Proposition}
\newtheorem{corollary}[theorem]{Corollary}
\theoremstyle{definition}
\newtheorem{definition}[theorem]{Definition}

\theoremstyle{remark}

\numberwithin{equation}{section}

\begin{document}
\setcounter{page}{1}
\title[Hereditary properties of character injectivity]
{Hereditary properties of character injectivity with applications to
semigroup algebras}

  \author[M. Essmaili, M. Fozouni and J. Laali]{M. Essmaili$^{*}$, M. Fozouni and J. Laali}

    \address{Faculty of Mathematical and Computer Science,
    Kharazmi University, 50 Taleghani Avenue, 15618 Tehran,
    Iran.}
    \email{m.essmaili@khu.ac.ir}

    \email{fozouni@khu.ac.ir}

    \email{laali@khu.ac.ir}



\begin{abstract}
In this paper, we investigate the notion $\phi$-injectivity for
Banach $A$-modules, where $\phi$ is a character on $A.$ We obtain
some hereditary properties of $\phi$-injectivity for certain classes
of Banach modules related to closed ideals. These results allow us
to study $\phi$-injectivity of certain Banach $A$-modules in
commutative case, specially  $\ell^{1}$-semilattice algebras. As an
application, we give an example of a non-injective Banach module
which is $\phi$-injective for each character $\phi.$
\end{abstract} \maketitle
\noindent
\section{introduction}
Suppose that $A$ is a Banach algebra. We denote by $\textbf{A-mod}$
and $\textbf{mod-A}$ the categories of Banach left $A$-modules and
Banach right $A$-modules, respectively. In the case where $A$ is
unital, we also denote by $\textbf{A-unmod}$ the categories of
unital Banach left $A$-modules. For each $E,F\in\textbf{A-mod},$ let
$_{A}B(E,F)$ be the closed subspace of $B(E,F)$ consisting of the
left $A$-module morphisms. An operator $T\in B(E,F)$ is called {\it
admissible} if $\text{ker} T$ and $\text{Im} T$ are closed
complemented subspaces of $E$ and $F$, respectively. It is easy to
verify that $T$ is admissible if and only if there exists $S\in
B(F,E)$ such that $T\circ S\circ T=T.$

A Banach left $A$-module $E$ is called {\it injective} if for each
$F, K\in\textbf{A-mod}$ and admissible monomorphism
$T\in_{A}\hspace{-0.1cm}B(F,K),$ the induced map
$_{A}B(K,E)\longrightarrow\hspace{-0.1cm}_{A}B(F,E)$ is onto. We
also say $E\in\textbf{mod-A}$ is {\it flat} if the dual module of
$E^{*}\in\textbf{A-mod}$ is injective with the following left module
action:
$$(a\cdot f)(x)=f(x\cdot a)\qquad (a\in A, x\in E).$$
The notions of injectivity and flatness of Banach algebras were
introduced by A. Ya. Helemskii. These notions have been studied for
various classes of Banach modules; see \cite{Dales. Polyakov},
\cite{Helem}, \cite{Ramsden Paper} and \cite{white} for more
details. Recently, Ramsden in \cite{Ramsden Paper} studied
injectivity and flatness of Banach modules over semigroup algebras.
It is well known that if $A$ is amenable, then every Banach
$A$-modules is flat but the converse is a long standing open
problem. We recall that the answer is positive for some classes of
Banach algebras associated with locally compact groups such as, the
class of group algebras and measure algebras; see \cite{Dales.
Polyakov} and \cite{Ramsden}.

Kaniuth, Lau and Pym introduced and studied in \cite{klp} and
\cite{klp2} the notion of $\phi$-amenability for Banach algebras,
where $\phi:A\longrightarrow\mathbb{C}$ is a character, i.e., a
non-zero homomorphism on $A$. Afterwards, Monfared introduced and
studied in \cite{m} the notion of character amenability for Banach
algebras. Let $\Delta(A)$ be the set of all characters of the Banach
algebra $A$, and let $\phi\in\Delta(A)$. The Banach algebra $A$ is
called {\it left $\phi$-amenable} if for all Banach $A$-bimodules
$E$ for which the left module action is given by
$$a\cdot x=\phi(a)x \qquad (a\in A, x\in E),$$ every
derivation $D:A\longrightarrow E^{*}$ is inner. It is clear that
amenability of $A$ implies $\phi$-amenability for all
$\phi\in\Delta(A).$

Recently, Nasr-Isfahani and Soltani Renani in \cite{Nasr} introduced
and studied the notion of $\phi$-injectivity and $\phi$-flatness for
Banach modules (see Definition $2.1$). As an important result, it is
shown in \cite[Proposition 3.1]{Nasr} that the Banach algebra $A$ is
left $\phi$-amenable if and only if every Banach left $A$-modules
$E$ is $\phi$-flat. Indeed, this result gives a positive answer to
the above open problem arises by A. Ya. Helemskii in this homology
setting based on character $\phi$. Furthermore, they obtained some
necessary and sufficient conditions for $\phi$-injectivity and
characterized $\phi$-injectivity of Banach modules in terms of a
coretraction problem; see \cite[Theorem 2.4]{Nasr}.

This paper is organized as follows. In Section 2, after recalling
some definitions, we investigate some properties of
$\phi$-injectivity for Banach modules. Indeed, we obtain a
sufficient condition for $\phi$-injectivity of Banach left
$A$-modules, in the case where $A$ is a commutative Banach algebra.
Moreover, we give some hereditary properties of $\phi$-injectivity
for Banach $A$-modules related to the closed ideals of Banach
algebra $A$. As the main result, we show that if $J$ is a left
invariant complemented ideal in $A$, then $\phi$-injectivity of $J$
and $A/J$ in $\textbf{A-mod}$ is equivalent to the
$\phi$-injectivity of $A$ in $\textbf{A-mod}$ (Theorem \ref{Cor:
Directsum}). In Section 3, by using the results of Section 2, we
study $\phi$-injectivity of certain $\ell^{1}$-semilattice algebras
and show that $\ell^{1}(\mathbb{N}_{\wedge})$ as a Banach left
$\ell^{1}(\mathbb{N}_{\wedge})$-module is $\phi$-injective for each
character $\phi,$ although is not injective.

\section{$\phi$-injectivity and some hereditary properties}
First, we recall some standard notations that we shall use and
define the notions of $\phi$-injectivity and $\phi$-flatness of
Banach modules.

Let $A$ be a Banach algebra and $E\in\textbf{A-mod}.$ Throughout the
paper, we regard $E$ as a Banach left $A^{\sharp}$-module (the
unitization of $A$) with the following left module action:
$$(a,\lambda)\cdot x=a\cdot x+\lambda x \qquad(a\in A, \lambda\in\mathbb{C}, x\in E).$$
Moreover, the space $B(A,E)$ is a Banach $A$-bimodule with the
following module actions:
$$(a\cdot T)(b)=T(ba), \qquad (T\cdot a)(b)=T(ab)\quad (T\in B(A,E), a, b\in A).$$
Suppose that $A$ is a Banach algebra and $\phi\in\Delta(A)$. For
each $E\in\textbf{A-mod}$ we define,
$$I(\phi,E)=\textrm{span}\{a\cdot x-\phi(a)x : a\in A, x\in E\}.$$
Following \cite{Nasr}, we also consider
$$_{\phi}B(A^{\sharp},E)=\{T\in B(A^{\sharp},E) :
T(ab-\phi(b)a)=a\cdot T(b-\phi(b)e^{\sharp})\:\:\text{for all}\:\:
a, b\in A\},$$ where $e^{\sharp}=(0,1)$ denotes the unite of
$A^{\sharp}.$ It is straightforward to check that
$_{\phi}B(A^{\sharp},E)$ is a closed $A$-submodule of
$B(A^{\sharp},E).$ Moreover, we define {\it the canonical morphism}
$_{\phi}\Pi^{\sharp}:E\longrightarrow _{\phi}B(A^{\sharp},E)$ as
follows:
$$_{\phi}\Pi^{\sharp}(x)(a)=a\cdot x\qquad (x\in E, a\in A^{\sharp}).$$
\begin{definition}
Let $A$ be a Banach algebra, $\phi\in \Delta(A)$ and $E\in
\textbf{A-mod}$. We say that $E$ is {\it $\phi$-injective} if, for
each $F, K\in \textbf{A-mod}$ and admissible monomorphism
$T:F\longrightarrow K$ with $I(\phi,K)\subseteq\text{Im}(T)$, the
induced map $T_{E}:\hspace{-0.1cm}_{A}B(K,E)\longrightarrow\
_{A}B(F,E)$ defined by $T_{E}(R)=R\circ T$ is onto.
\end{definition}
The following theorem gives a characterization of $\phi$-injectivity
in terms of a coretraction problem.
\begin{theorem}$($\cite[Theorem 2.4]{Nasr}$)$
Let $A$ be a Banach algebra and $\phi\in\Delta(A)$. For
$E\in\textbf{A-mod}$ the following statements are equivalent.
\begin{enumerate}
\item[(i)] $E$ is $\phi$-injective.
\item[(ii)] $_{\phi}\Pi^{\sharp}\in\hspace{-0.1cm}_{A}B(E,_{\phi}B(A^{\sharp},E))$ is a coretraction, $($that is there exists
$_{\phi}\rho^{\sharp}\in\: _{A}B(_{\phi}B(A^{\sharp},E), E)$ such
that is a left inverse for $_{\phi}\Pi^{\sharp})$.
\end{enumerate}
\end{theorem}
A Banach right (left) $A$-module $E$ is $\phi$-flat if $E^{*}$ is
$\phi$-injective as a left (right) $A$-module. It is shown that
Banach algebra $A$ is left $\phi$-amenable if and only if each
Banach left $A$-module $E$ is $\phi$-flat \cite[Proposition
3.1]{Nasr}.

In this section, we give some hereditary properties of
$\phi$-injectivity for certain classes of Banach modules. We also
consider some hereditary properties of $\phi$-injectivity of Banach
left $A$-modules with their ideals. We first give a sufficient
condition for $\phi$-injectivity of Banach left $A$-module $E$ in
the case where $A$ is a commutative Banach algebra. Following
\cite[Definition 1.4.4]{Dales}, {\it the annihilator} of $E$ is
defined by $E^{\perp}=\{a\in A: a\cdot E=\{0\}\}$.
\begin{theorem}\label{Th: 2}
Let $A$ be a commutative Banach algebra, $\phi\in \Delta(A)$ and
$E\in \textbf{A-mod}.$ If
$E^{\perp}\cap(\ker(\phi))^{c}\neq\emptyset$, then $E$ is
$\phi$-injective.
\end{theorem}
\begin{proof}
Let $a_{0}\in E^{\perp}\cap(\ker(\phi))^{c}$. We can assume that
$\phi(a_{0})=1$ and define the map
$_{\phi}\rho^{\sharp}:$$_{\phi}B(A^{\sharp},E)\longrightarrow E$ by
\begin{center}
$_{\phi}\rho^{\sharp}(T)=T(e^{\sharp}-a_{0})\hspace{0.5cm}(T\in$$
_{\phi}B(A^{\sharp},E)).$
\end{center}
Hence, for each $x\in E$ we have
\begin{center}
$_{\phi}\rho^{\sharp}\circ$$_{\phi}\Pi^{\sharp}(x)=$$_{\phi}\Pi^{\sharp}(x)(e^{\sharp}-a_{0})=(e^{\sharp}-a_{0})\cdot
x=x.$
\end{center}
Therefore, $_{\phi}\rho^{\sharp}\circ$$_{\phi}\Pi^{\sharp}=I_{E}$.
On the other hand, for each $a\in A$ and $T\in$$
_{\phi}B(A^{\sharp},E)$ we have
\begin{equation}
\begin{split}
_{\phi}\rho^{\sharp}(a\cdot T)&=(a\cdot T)(e^{\sharp}-a_{0})=T((e^{\sharp}-a_{0})\cdot a)\\
&=T(a-a_{0}a)\\
&=T(a\phi(a_{0})-aa_{0}).
\end{split}
\end{equation}
Since $T\in$$ _{\phi}B(A^{\sharp},E)$ we have
$T(a\phi(a_{0})-aa_{0})=a\cdot T(e^{\sharp}-a_{0}).$ Now, using
($2.1$) we conclude that
\begin{align*}
_{\phi}\rho^{\sharp}(a\cdot T)=a\cdot
T(e^{\sharp}-a_{0})=a\cdot_{\phi}\rho^{\sharp}(T).
\end{align*}
It follows that $_{\phi}\rho^{\sharp}$ is a left $A$-module
morphism. Hence, $E$ is a $\phi$-injective Banach left $A$-module.
\end{proof}
\begin{corollary}\label{Cor: 2} Let $A$ be a commutative Banach algebra, and $J$ be a closed ideal of $A$
such that $\phi_{|J}\neq 0$. Then $A/J$ is $\phi$-injective as a
Banach left $A$-module.
\end{corollary}
\begin{proof}
Since $\phi_{|J}\neq 0$, it is easy to check that
$(A/J)^{\perp}\cap(\ker(\phi))^{c}\neq \emptyset$. Now, apply
Theorem \ref{Th: 2}.
\end{proof}
\begin{corollary}\label{Th: Inj Module}
Let $A$ be a commutative Banach algebra, $\phi\in\Delta(A)$ and let $E\in\textbf{A-mod}$ with $I(\phi,E)=\{0\}$. Then for all $\psi\in\Delta(A)\setminus\{\phi\}$, $E\in \textbf{A-mod}$ is $\psi$-injective.
\end{corollary}
\begin{proof} Since $\phi\neq\psi$ there exists $a_{0}\in A$ such that $\phi(a_{0})=0$ and $\psi(a_{0})=1$.
On the other hand, since $I(\phi,E)=\{0\}$ we conclude that $a_{0}\in E^{\perp}\cap(\ker(\psi))^{c}$ and the proof
is complete.
\end{proof}
Now, we give some hereditary properties of $\phi$-injectivity of
Banach modules that we shall use. Recall that $E\in\textbf{A-mod}$
is {\it faithful} if $ A\cdot x=\{0\}$ implies that $x=0.$
\begin{theorem}\label{Th: Her,of,Id}
Let $A$ be a Banach algebra,$E\in\textbf{A-mod}$, $\phi\in\Delta(A)$ and $J$ be a closed
ideal of $A$ such that $\phi_{|J}\neq 0$.
\begin{enumerate}
\item[(i)] Suppose that $J$ has an identity and $E\in\textbf{J-unmod}$. If $E\in\textbf{A-mod}$ is
$\phi$-injective, then $E\in\textbf{J-unmod}$ is
$\phi_{|J}$-injective.
\item[(ii)] If $E\in\textbf{J-mod}$ is $\phi_{|J}$-injective and
faithful, then $E\in\textbf{A-mod}$ is $\phi$-injective.
\end{enumerate}
\end{theorem}
\begin{proof}(i) Suppose that $E\in \textbf{A-mod}$ is $\phi$-injective. Let $F$ and $K$ be in $\textbf{J-mod}$ and $T:F\longrightarrow K$ be an
admissible monomorphism with $I(\phi_{|J},K)\subseteq \textrm{Im}T$.
We claim that the induced map, $_{J}B(K,E)\longrightarrow\
_{J}B(F,E)$ defined by, $R\longrightarrow R\circ T$ is onto. Suppose
that $e_{J}$ is the identity of $J$. We can consider $F$ and $K$ as
Banach left $A$-modules with the following module actions:
\begin{align*}
&a\bullet f=(ae_{J})\cdot f\hspace{0.5cm}(a\in A, f\in F),\\
&a\bullet k=(ae_{J})\cdot k\hspace{0.5cm}(a\in A, k\in K).
\end{align*}
We denote these $A$-modules with $\widetilde{F}$ and
$\widetilde{K},$ respectively. Take $W\in$$_{J}B(F,E)$ and define
the map $\widetilde{W}:\widetilde{F}\longrightarrow E$ by
$\widetilde{W}(f)=W(f)$. For each $a\in A$ and $f\in F$ we have,
\begin{align*}
\widetilde{W}(a\bullet f)&=W((ae_{J})\cdot f)=(ae_{J})\cdot
W(f)\\
&=a\cdot (e_{J}\cdot W(f))=a\cdot W(f)\\
&=a\cdot\widetilde{W}(f).
\end{align*}
So $\widetilde{W}$ is a left $A$-module morphism. Moreover, the map
$\widetilde{T}:\widetilde{F}\longrightarrow\widetilde{K}$ defined by
$\widetilde{T}(f)=T(f)$ is an admissible monomorphism such that
\begin{align*}
I(\phi,\widetilde{K})&=\textrm{span}\{a\bullet k-\phi(a)k : a\in A, k\in K\}\\
&=\textrm{span}\{(ae_{J})\cdot k-\phi(ae_{J})k : a\in A, k\in K\}\\
&\subseteq \textrm{Im}T=\textrm{Im}\widetilde{T}.
\end{align*}
Since $E\in \textbf{A-mod}$ is $\phi$-injective, there exist
$S\in$$_{A}B(\widetilde{K},E)$ such that $S\circ
\widetilde{T}=\widetilde{W}$. On the other hand, for each $a\in J$
and $k\in K$ we have
\begin{equation*}
a\cdot S(k)=S(a\bullet k)=S((ae_{J})\cdot k)=S(a\cdot k).
\end{equation*}
It follows that $S\in$$_{J}B(K,E)$. Now, we conclude that $E\in
\textbf{J-unmod}$ is $\phi_{|J}$-injective.

(ii) Let $F$ and $K$ be in $\textbf{A-mod}$ and $T:F\longrightarrow
K$ be an admissible monomorphism and take
$W\in\hspace{-0.1cm}_{A}B(F,E)$. So $W\in\hspace{-0.1cm}_{J}B(F,E)$
and there exists $S\in\hspace{-0.1cm}_{J}B(K,E)$ such that $S\circ
T=W$. For each $a\in J$, $b\in A$ and $k\in K$, we have
\begin{align*}
a\cdot(S(b\cdot k)-b\cdot S(k))&=a\cdot S(b\cdot k)-(ab)\cdot S(k)\\
&=S(ab\cdot k)-S(ab\cdot k)=0.
\end{align*}
Since $E\in \textbf{J-mod}$ is faithful, we conclude that $S(b\cdot
k)=b\cdot S(k).$ It follows that $S\in\hspace{-0.1cm}_{A}B(K,E)$ and
the proof is complete.
\end{proof}
\begin{corollary}
Let $A$ be a Banach algebra, $\phi\in\Delta(A)$ and $J$ be a closed
ideal of $A$ with an identity such that $\phi_{|J}\neq 0$. Then
$J\in\textbf{A-mod}$ is $\phi$-injective if and only if
$J\in\textbf{J-mod}$ is $\phi_{|J}$-injective.
\end{corollary}
B. E. Forrest in \cite{Forrest} introduced the notion of invariantly
complemented submodules in categories of Banach modules. In the
sequel, we obtain some results for $\phi$-injectivity of invariantly
complemented ideals.
\begin{definition}(\cite[Definition
6.3]{Forrest}) Let $X$ be a Banach left $A$-module and $Y$ be a
Banach $A$-submodule of $X$. We say that $Y$ is {\it left
$($right$)$ invariantly complemented} in $X$ if there exists
$P\in\hspace{-0.1cm}_{A}B(X,Y)$ ($P\in B_{A}(X,Y)$) such that
$P^{2}=P$ and $P(X)=Y$.
\end{definition}
\begin{theorem}\label{Th: 3} Let $\{E_{\alpha}\}_{\alpha\in \Gamma}$ be a collection of Banach left $A$-modules and consider
$E=\ell^{1}-\bigoplus_{\alpha\in \Gamma} E_{\alpha}$ as a Banach
left $A$-module with the natural module action.
\begin{enumerate}
\item[(i)] If $E$ is $\phi$-injective, then for each $\alpha\in\Gamma,$ $E_{\alpha}$ is $\phi$-injective.
\item[(ii)] Conversely, if $\Gamma$ is finite and each $E_{\alpha}$ is $\phi$-injective, then $E$ is $\phi$-injective.
\end{enumerate}
\end{theorem}
\begin{proof} (i) It is obvious that each $E_{\alpha}$ is left invariantly complemented in $E$. Hence, for each $\alpha\in\Gamma$,
let $P_{\alpha}\in$$_{A}B(E,E_{\alpha})$ such that
$P_{\alpha}(E)=E_{\alpha}$ and $P^{2}_{\alpha}=P_{\alpha}$. Also,
let $i_{\alpha}:E_{\alpha}\longrightarrow E$ be the natural
embedding of $E_{\alpha}$ into $E$.

Let $E$ be $\phi$-injective. Then there exists
$_{\phi}\rho^{E}\in$$_{A}B(_{\phi}B(A^{\sharp},E),E)$ such that is a
left inverse for
$_{\phi}\Pi^{E}:E\longrightarrow\hspace{-0.1cm}_{\phi}B(A^{\sharp},E)$.
For each $\alpha\in \Gamma$, we define the map
$_{\phi}\rho^{\alpha}:$$_{\phi}B(A^{\sharp},E_{\alpha})\longrightarrow
E_{\alpha}$ by
\begin{center}
$_{\phi}\rho^{\alpha}(T)=P_{\alpha}\circ$$_{\phi}\rho^{E}(i_{\alpha}\circ
T)\qquad(T\in$$_{\phi}B(A^{\sharp},E_{\alpha})).$
\end{center}
We claim that $_{\phi}\rho^{\alpha}$ is a left $A$-module morphism
such that
$_{\phi}\rho^{\alpha}\circ$$_{\phi}\Pi^{\alpha}=I_{E_{\alpha}}$.
Indeed, since for each $x\in E_{\alpha}$,
$i_{\alpha}\circ$$_{\phi}\Pi^{\alpha}(x)=$$_{\phi}\Pi^{E}(i_{\alpha}(x))$,
so we have
\begin{align*}
_{\phi}\rho^{\alpha}\circ_{\phi}\Pi^{\alpha}(x)&=P_{\alpha}\circ_{\phi}\rho^{E}(i_{\alpha}\circ
_{\phi}\Pi^{\alpha}(x))=P_{\alpha}\circ_{\phi}\rho^{E}(_{\phi}\Pi^{E}(i_{\alpha}(x)))\\
&=P_{\alpha}(i_{\alpha}(x))=x.
\end{align*}
Therefore,
$_{\phi}\rho^{\alpha}\circ$$_{\phi}\Pi^{\alpha}=I_{E_{\alpha}}$. On
the other hand, since $P_{\alpha}\in$$_{A}B(E,E_{\alpha})$ it is
easy to see that $_{\phi}\rho^{\alpha}$ is a left $A$-module
morphism and the proof is complete.

(ii) Suppose that for each $\alpha\in \Gamma$, $E_{\alpha}$  is
$\phi$-injective. So, for each $\alpha\in \Gamma$ there exists
$_{\phi}\rho^{\alpha}\in$$_{A}B(_{\phi}B(A^{\sharp},E_{\alpha}),E_{\alpha})$
for which
$_{\phi}\rho^{\alpha}\circ$$_{\phi}\Pi^{\alpha}=I_{E_{\alpha}}$.
Define the map $\rho:$$_{\phi}B(A^{\sharp},E)\longrightarrow E$ by
\begin{center}
$\rho(T)=(_{\phi}\rho^{\alpha}(P_{\alpha}\circ T))_{\alpha\in
\Gamma}\qquad(T\in$$_{\phi}B(A^{\sharp},E)).$
\end{center}
Since $\Gamma$ is finite, $\rho$ is well-defined. For each $a\in A$ and $T\in$$_{\phi}B(A^{\sharp},E)$ we have
\begin{align*}
\rho(a\cdot T)&=(_{\phi}\rho^{\alpha}(P_{\alpha}\circ (a\cdot T)))_{\alpha\in \Gamma}=(_{\phi}\rho^{\alpha}(a\cdot(P_{\alpha}\circ T)))_{\alpha\in \Gamma}\\
&=(a\cdot_{\phi}\rho^{\alpha}(P_{\alpha}\circ T))_{\alpha\in \Gamma}=a\cdot(_{\phi}\rho^{\alpha}(P_{\alpha}\circ T))_{\alpha\in \Gamma}\\
&=a\cdot\rho(T).
\end{align*}
Moreover, if $x=(x_{\alpha})_{\alpha\in\Gamma}$ is an arbitrary
element of $E$, it is easy to see that
$P_{\alpha}\circ$$_{\phi}\Pi^{E}(x)=$$_{\phi}\Pi^{\alpha}(x_{\alpha})$.
Hence,
\begin{align*}
\rho\circ_{\phi}\Pi^{E}(x)&=(_{\phi}\rho^{\alpha}(P_{\alpha}\circ
_{\phi}\Pi^{E}(x)))_{\alpha\in
\Gamma}=(_{\phi}\rho^{\alpha}(_{\phi}\Pi^{\alpha}(x_{\alpha})))_{\alpha\in
\Gamma}\\
&=(x_{\alpha})_{\alpha\in\Gamma}=x.
\end{align*}
Therefore, we conclude that $E$ is $\phi$-injective.
\end{proof}
\begin{theorem}\label{Cor: Directsum}
Let $A$ be a Banach algebra, $\phi\in \Delta(A)$, $B$ be a
subalgebra of $A$ and $J$ be a closed left ideal of $A$. Then the
following assertions holds:
\begin{enumerate}
\item[(i)] If $B$ is left invariantly complemented in $A$ and $A$  is $\phi$-injective in
\textbf{A-mod},  then $B$ is $\phi$-injective in \textbf{A-mod}.
\item[(ii)] If $J$ is left invariantly complemented, then $J$ and $A/J$ are $\phi$-injective in \textbf{A-mod}
if and only if $A$ is $\phi$-injective in \textbf{A-mod}.
\end{enumerate}
\end{theorem}
\begin{proof}(i) Since $B$ is left invariantly complemented in $A$,
there exists an onto projection $P\in$$_{A}B(A,B)$. Hence
$$A\cong\textrm{Im} P\oplus\ker P=B\oplus\ker P,$$ as a Banach
left $A$-module. Therefore, by Theorem \ref{Th: 3} it follows that
$B$ is $\phi$-injective in $\textbf{A-mod}.$

(ii) Since $J$ is a left invariant complemented ideal in $A,$ there
exists an onto projection $P\in$$_{A}B(A,J)$. We claim that $A\cong
J\oplus \frac{A}{J}$ as a Banach left $A$-module. To see this,
define the map $T:A\longrightarrow J\oplus \frac{A}{J}$ by
\begin{equation*}
T(a)=(P(a), a+J)\qquad(a\in A).
\end{equation*}
First, $T$ is a left $A$-module morphism, because for each $a, b\in
A$ we have,
\begin{align*}
T(ab)&=(P(ab),ab+J)=(aP(b),a(b+J))\\
&=a\cdot (P(b),b+J)=a\cdot T(b).
\end{align*}
On the other hand, if $a\in \textrm{Im}P\ \cap\ \ker P$, then there
exists $b\in A$ such that $P(b)=a$. Hence, $a=P(b)=P(P(b))=P(a)=0$.
This follows that $\textrm{Im}P\cap \ker P=\{0\}$ and so $T$ is
one-to-one. Moreover, $T$ is onto because for each $(a,b+J)\in
J\oplus \frac{A}{J}$ if we put $c=a+b-P(b)$, then $T(c)=(a,b+J)$.
Now, the result follows from Theorem \ref{Th: 3}.
\end{proof}
As an application of Theorem \ref{Cor: Directsum} and Corollary
\ref{Cor: 2}, we have the following result for commutative Banach
algebras.
\begin{corollary}\label{Cor: A and I. A is Com}
Let $A$ be a commutative Banach algebra, $\phi\in \Delta(A)$ and let
$J$ be a closed invariant complemented ideal of $A$ such that
$\phi_{|J}\neq 0$. Then $A\in\textbf{A-mod}$ is $\phi$-injective if
and only if $J\in\textbf{A-mod}$ is $\phi$-injective.
\end{corollary}
\section{Applications to semigroup algebras}
In this section, we apply our later results to study
$\phi$-injectivity of certain commutative semigroup algebras and
give some examples of non-injective modules which are
$\phi$-injective for each character $\phi.$ First, we need some
basic facts about semigroup algebras.

Let $S$ be a semigroup and let $E(S)=\{s\in S : s^{2}=s\}.$ We say
that $S$ is a {\it semilattice} if $S$ is commutative and $E(S)=S.$
A {\it semi-character} on $S$ is a non-zero homomorphism
$\widehat{\phi}:S\longrightarrow\{z\in\mathbb{C}: |z|\leq 1\}.$ The
space of semi-characters on $S$ is denoted by $\Phi_{S}.$ The
semi-character $\widehat{\phi}_{S}:S\longrightarrow\{z\in\mathbb{C}:
|z|\leq 1\}$, defined by
$$\widehat{\phi}_{S}(t)=1 \qquad (t\in S),$$ is called the {\it
augmentation character} on $S$. For each
$\widehat{\phi}\in\Phi_{S},$ we associate the map
$\phi:\ell^{1}(S)\longrightarrow \mathbb{C}$ defined by
$$\phi(f)=\sum_{s\in S}\widehat{\phi}(s)f(s) \qquad (f\in\ell^{1}(S)).$$
It is easy to verify that $\phi\in\Delta(\ell^{1}(S))$ and every
character on $\ell^{1}(S)$ arises in this way. Indeed, we have
$$\Delta(\ell^{1}(S))=\{\phi: \widehat{\phi}\in\Phi_{S}\}.$$
We also define the convolution of two elements $f,g\in \ell^{1}(S)$
by
$$(f\ast g)(s)= \sum_{uv=s}f(u)g(v) \qquad (s\in S),$$
where $\sum_{uv=s}f(u)g(v)=0,$ when there are no elements $u,v\in S$
with $uv=s.$ Then $(\ell^{1}(S),\ast,{\|\cdot\|}_{1})$ becomes a
Banach algebra that is called the {\it semigroup algebra} of $S.$
The following proposition immediately follows from Corollary
\ref{Cor: A and I. A is Com}.
\begin{proposition}\label{pp}
Let $S$ be a semilattice, $\phi\in\Delta(\ell^{1}(S))$ and $I$ be a
closed invariant complemented ideal in $\ell^{1}(S)$ such that
$\phi_{|_{I}}\neq 0.$ Then $\ell^{1}(S)\in\ell^{1}(S)\textbf{-mod}$
is $\phi$-injective if and only if $I\in\ell^{1}(S)\textbf{-mod}$ is
$\phi$-injective.
\end{proposition}
Let $\ell^{1}(\mathbb{N}_{\wedge})$ be the semigroup algebra on semigroup $S=(\mathbb{N},{\wedge})$ with the following
product:
\begin{equation*}
\mathbb{N}\times \mathbb{N}\longrightarrow \mathbb{N},\quad
(m,n)\longrightarrow m\wedge n=\min\{m, n\}.
\end{equation*}
It is easy to check that $\Phi_{S}=\{\widehat{\phi}_{n} :
n\in\mathbb{N}\},$ where for each $n\in\mathbb{N},$
$$\begin{array}{lll}
\widehat{\phi}_{n}(m) =\left\{\begin{array}{l} 1 \quad\text{if}\quad
m\geq n
\\ 0 \quad \text{if} \quad m< n
\end{array}\qquad (m\in\mathbb{N}).\right.
\end{array}$$
For each $n\in \mathbb{N}$, let $I_{n}$ be the ideal of
$\ell^{1}(\mathbb{N}_{\wedge})$ generated by the set $\{\delta_{1},
\delta_{2}, \delta_{3},\ldots, \delta_{n}\}$. It is easy to see that
$\ell^{1}(\mathbb{N}_{\wedge})/I_{n}$ does not have an identity. So
$I_{n}$ does not have a modular identity and using \cite[Corollary
2.2.8 (ii)]{Ramsden}, we conclude that
$\ell^{1}(\mathbb{N}_{\wedge})/I_{n}$ is not injective as a Banach
left $\ell^{1}(\mathbb{N}_{\wedge})$-module. Furthermore, we recall
that $\ell^{1}(\mathbb{N}_{\wedge})$ as a Banach left
$\ell^{1}(\mathbb{N}_{\wedge})$-module is not injective, because it
does not have a right identity \cite[Example 4.10]{Dales. Lau}. As
mentioned above, we regard the map
$\phi_{n}:\ell^{1}(\mathbb{N}_{\wedge})\longrightarrow \mathbb{C}$
as a character on $\ell^{1}(\mathbb{N}_{\wedge})$ which is defined
by
\begin{equation*}
\phi_{n}(f)=\sum_{i=n}^{\infty}f(i)\qquad(f\in
\ell^{1}(\mathbb{N}_{\wedge})).
\end{equation*}
The following theorem shows that $\ell^{1}(\mathbb{N}_{\wedge})$ as
a left $\ell^{1}(\mathbb{N}_{\wedge})$-module is $\phi$-injective
for each $\phi\in \Delta(\ell^{1}(\mathbb{N}_{\wedge}))$, although
is not injective.
\begin{theorem}\label{Th: L1Min is 1-inj}
With the above notations, we have following assertions:
\begin{enumerate}
\item[(i)] $\ell^{1}(\mathbb{N}_{\wedge})/I_{n}$ as a Banach left
$\ell^{1}(\mathbb{N}_{\wedge})$-module is $\phi_{n}$-injective, for
each $n\in\mathbb{N}$.
\item[(ii)] $\ell^{1}(\mathbb{N}_{\wedge})$ as a Banach left
$\ell^{1}(\mathbb{N}_{\wedge})$-module is $\phi_{n}$-injective,  for
each $n\in\mathbb{N}$.
\end{enumerate}
\end{theorem}
\begin{proof}(i) Since $\ell^{1}(\mathbb{N}_{\wedge})$ is commutative and
$(\phi_{n})_{|I_{n}}\neq 0$, by Corollary \ref{Cor: 2}, it follows
that $\ell^{1}(\mathbb{N}_{\wedge})/I_{n}$ as a Banach left
$\ell^{1}(\mathbb{N}_{\wedge})$-module is $\phi_{n}$-injective.

(ii) First, we show that for each $n\in \mathbb{N}$, $I_{n}$ is an
invariant complemented ideal of $\ell^{1}(\mathbb{N}_{\wedge})$. To
see this, suppose that the map
$P_{n}:\ell^{1}(\mathbb{N}_{\wedge})\longrightarrow I_{n}$ is
defined by
\begin{equation*}
P_{n}(f)=\sum_{i=1}^{n-1}f(i)\delta_{i}+(\sum_{i=n}^{\infty}f(i))\delta_{n}\qquad(f\in
\ell^{1}(\mathbb{N}_{\wedge})).
\end{equation*}
It is easy to check that $P_{n}$ is a projection on $I_{n}.$
Moreover, if $f$ or $g$ belong to $I_{n}$ we have $P_{n}(f\ast
g)=f\ast P_{n}(g)$. Now suppose that $f, g$ are not in $I_{n}$. We
can suppose that $f=\delta_{i}$ and $g=\delta_{j}$ such that
$n<i\leq j$. Hence, we have
\begin{align*}
P_{n}(\delta_{i}\ast\delta_{j})=P_{n}(\delta_{i})=\delta_{n}=\delta_{i}\ast\delta_{n}=\delta_{i}\ast
P_{n}(\delta_{j}).
\end{align*}
This follows that $I_{n}$ is an invariant complemented closed ideal
of $\ell^{1}(\mathbb{N}_{\wedge})$.

$\mathbf{Case\ 1}$: We prove that $A=\ell^{1}(\mathbb{N}_{\wedge})$
is $\phi_{1}$-injective. Since $I_{1}$ is invariantly complemented
in $A,$ $(\phi_{1})_{|I_{1}}\neq 0$ and $A$ is commutative, by
Proposition \ref{pp}, it suffices to show that
$I_{1}\in\textbf{A-mod}$ is $\phi_{1}$-injective.

Define the map $\rho:\ _{\phi_{1}}B(A^{\sharp},I_{1})\longrightarrow
I_{1}$ by
$$\rho(T)=T(e^{\sharp})\qquad  (T\in\
_{\phi_{1}}B(A^{\sharp},I_{1}).$$
Clearly $\rho$ is a left inverse
for $\ _{\phi_{1}}\Pi^{\sharp}:I_{1}\longrightarrow \
_{\phi_{1}}B(A^{\sharp},I_{1})$. We claim that $\rho$ is a left
$A$-module morphism.

For each $f\in A$ and $g\in I_{1}$ we note that $f\ast
g=\phi_{1}(f)g$. Hence, if $T\in\ _{\phi_{1}}B(A^{\sharp},I_{1})$,
then for each $f, g\in A$ we have
\begin{equation*}
T(f\ast g-\phi_{1}(g)f)=f\cdot T(g-\phi_{1}(g)e^{\sharp}).
\end{equation*}
In the case where $g=\delta_{1}$, we conclude that
\begin{align*}
T(\phi_{1}(f)\delta_{1}-\phi_{1}(g)f)&=T(f\ast g-\phi_{1}(g)f)\\
&=f\cdot T(g-\phi_{1}(g)e^{\sharp})\\
&=f\ast T(g-\phi_{1}(g)e^{\sharp})\\
&=\phi_{1}(f)T(\delta_{1}-\phi_{1}(\delta_{1})e^{\sharp})\\
&=\phi_{1}(f)T(\delta_{1}-e^{\sharp}).
\end{align*}
Therefore, $T(f)=\phi_{1}(f)T(e^{\sharp})$. This follows that for
each $f\in A$ and $T\in\ _{\phi_{1}}B(A^{\sharp},I_{1})$ we have
\begin{align*}
\rho(f\cdot T)&=(f\cdot T)(e^{\sharp})=T(e^{\sharp}f)=T(f)\\
&=\phi_{1}(f)T(e^{\sharp}) =f\cdot T(e^{\sharp})=f\cdot\rho(T).
\end{align*}
Hence, $\rho$ is a left $A$-module morphism which is also a left
inverse for $_{\phi_{1}}\Pi^{\sharp}$. Therefore,
$\ell^{1}(\mathbb{N}_{\wedge})$ is $\phi_{1}$-injective as a Banach
left $\ell^{1}(\mathbb{N}_{\wedge})$-module.

$\mathbf{Case\ 2}$: We prove that  $\ell^{1}(\mathbb{N}_{\wedge})$ is
$\phi_{n}$-injective as a Banach left
$\ell^{1}(\mathbb{N}_{\wedge})$-module for each  $n\geq 2$. Since
$I_{n}\in\textbf{I}_{n}\textbf{-mod}$ has the identity $\delta_{n}$,
so it is faithful. By Theorem \ref{Th: Her,of,Id}, it suffices to
show that $I_{n}\in\textbf{I}_{n}\textbf{-mod}$ is
$(\phi_{n})_{|I_{n}}$-injective. Let $\widehat{\phi_{n}}$ be the
semi-character on $\mathbb{N}_{\wedge}$ associated with $\phi_{n}$.
Since $S=(\{1,2,3,...,n\},\wedge)$ is a semilattice and
$\widehat{\phi_{n}}(1)=0$, it follows from \cite[Theorem
2.1]{Essmaili} that $I_{n}=\ell^{1}(\{1,2,3,\ldots,n\},\wedge)$ is
$(\phi_{n})_{|I_{n}}$-amenable. Therefore, we conclude that
$c_{0}(S)$ is $(\phi_{n})_{|I_{n}}$-flat in
$\textbf{mod-}\textbf{I}_{n}$ by \cite[Proposition 3.1]{Nasr}.
Hence, $I_{n}=c_{0}(S)^{*}$  is $(\phi_{n})_{|I_{n}}$-injective in
$\textbf{I}_{n}\textbf{-mod}$ and the proof is complete.
\end{proof}
Let $\ell^{1}(\mathbb{N}_{\vee})$ be the semigroup algebra on
semigroup $S=(\mathbb{N},{\vee})$ with the following product:
\begin{equation*}
\mathbb{N}\times \mathbb{N}\longrightarrow \mathbb{N},\quad
(m,n)\longrightarrow m\vee n=\max\{m, n\}.
\end{equation*}
It is easy to check that $\Phi_{S}=\{\widehat{\psi}_{n} :
n\in\mathbb{N}\}\cup\{\widehat{\psi_{S}}\},$ where
$\widehat{\psi_{S}}$ is the augmentation character and for each
$n\in\mathbb{N},$
$$\begin{array}{lll}
\widehat{\psi}_{n}(m) =\left\{\begin{array}{l} 1 \quad\text{if}\quad
m\leq n
\\ 0 \quad \text{if} \quad m> n
\end{array}\qquad (m\in\mathbb{N}).\right.
\end{array}$$
In \cite[Example 5.6]{Ramsden Paper}, it is proved that
$\ell^{1}(\mathbb{N}_{\vee})$ is not injective as a Banach left
$\ell^{1}(\mathbb{N}_{\vee})$-module. In the following theorem, we
show that $\ell^{1}(\mathbb{N}_{\vee})$ is $\phi$-injective for each
$\phi\in \Delta(\ell^{1}(\mathbb{N}_{\vee}))$.
\begin{theorem}\label{aa}
$\ell^{1}(\mathbb{N}_{\vee})$ as a Banach left
$\ell^{1}(\mathbb{N}_{\vee})$-module is $\phi$-injective for each
$\phi\in \Delta(\ell^{1}(\mathbb{N}_{\vee}))$.
\end{theorem}
\begin{proof} By \cite[Corollary 2.2]{Essmaili}, it follows that $\ell^{1}(\mathbb{N}_{\vee})$ is character amenable . Hence,
for each $\phi\in\Delta(\ell^{1}(\mathbb{N}_{\vee}))$,
$\ell^{1}(\mathbb{N}_{\vee})$ is $\phi$-amenable. On the other hand,
since $\mathbb{N}_{\vee}$ is weakly cancellative by \cite[Theorem
4.6]{Dales. Lau}, we conclude that $c_{0}(\mathbb{N}_{\vee})$ is a
Banach $\ell^{1}(\mathbb{N}_{\vee})$-module . Hence,
$c_{0}(\mathbb{N}_{\vee})$ is $\phi$-flat as a Banach right
$\ell^{1}(\mathbb{N}_{\vee})$-module \cite[Proposition 3.1]{Nasr}.
This follows that
$c_{0}(\mathbb{N}_{\vee})^{*}=\ell^{1}(\mathbb{N}_{\vee})$ is
$\phi$-injective as a Banach left
$\ell^{1}(\mathbb{N}_{\vee})$-module.
\end{proof}
For each $n\in \mathbb{N}$, let $J_{n}$ be the closed ideal of
$A=\ell^{1}(\mathbb{N}_{\vee})$ generated by
$\{\delta_{n},\delta_{n+1},\ldots\}$. It is easy to see that $J_{n}$
is invariantly complemented in $\ell^{1}(\mathbb{N}_{\vee})$.
Indeed, it is sufficient to consider the map $Q_{n}:A\longrightarrow
J_{n}$ defined by
$$Q_{n}(f)=(\sum_{i=1}^{n}f(i))\delta_{n}+\sum_{i=n+1}^{\infty}f(i)\delta_{i} \qquad(f\in A).$$
It is straightforward to check that $Q_{n}$ is an onto projection in
$_{A}B(A,J_{n})$.

As a consequence of Corollary \ref{Cor: A and I. A is Com} and
Theorem \ref{aa}, we give the following result.
\begin{corollary} For each $n\in \mathbb{N}$, $J_{n}$ as a Banach left $\ell^{1}(\mathbb{N}_{\vee})$-module is
$\phi$-injective for each
$\phi\in\Delta(\ell^{1}(\mathbb{N}_{\vee}))$ with $\phi_{|J_{n}}\neq
0$.
\end{corollary}
\bibliographystyle{amsplain}

\end{document}